\documentclass[a4paper, 12pt, reqno]{amsart}
\usepackage[matrix, arrow, curve]{xy}
\usepackage{amscd, euscript}

\newcommand{\xref}[1]{{\rm \ref{#1}}}

\newcommand{\PP}{\mathbb P}
\newcommand{\QQ}{\mathbb Q}

\newcommand{\CC}{\mathbb C}
\newcommand{\ZZ}{\mathbb Z}
\newcommand{\FF}{\mathbb F}
\newcommand{\OOO}{\mathcal O}

\newcommand{\gr}{\muu_2}
\newcommand{\HHH}{{\EuScript{H}}}
\newcommand{\LLL}{{\EuScript{L}}}
\newcommand{\EEE}{{\EuScript{E}}}

\newcommand{\AAA}{{\EuScript{A}}}
\newcommand{\FFF}{{\EuScript{F}}}

\newcommand{\down}[1]{\left\lfloor #1\right\rfloor}

\newcommand{\Mori}{\overline{\operatorname{NE}}}
\newcommand{\Pic}{\operatorname{Pic}}
\newcommand{\Pici}{\operatorname{Pic}^{\left\langle\tau\right\rangle}}
\newcommand{\Supp}{\operatorname{Supp}}
\newcommand{\Diff}{\operatorname{Diff}}
\newcommand{\Bs}{\operatorname{Bs}}
\newcommand{\muu}{\mbox{\boldmath $\mu$}}

\newtheorem{theorem}[subsection]{Theorem}
\newtheorem{proposition}[subsection]{Proposition}
\newtheorem{lemma}[subsection]{Lemma}
\newtheorem{claim}[subsection]{Claim}
\newtheorem{corollary}[subsection]{Corollary}
\theoremstyle{definition}

\newtheorem{remark}[subsection]{Remark}

\title{On Fano-Enriques threefolds}
\author{Yuri Prokhorov}
\thanks {Partially supported by grants CRDF-RUM, No. 1-2692-MO-05 and
RFBR, No. 05-01-00353-a.}
\address{Department of Algebra, Faculty of Mathematics, Moscow State
University, Moscow 117234, Russia}
\email{prokhoro@mech.math.msu.su}

\begin{document}
\begin{abstract}
Let $U\subset \mathbb P^N$ be a projective 
variety 
which is not a cone and 
whose hyperplane sections are smooth Enriques
surfaces.
We prove that the degree of such $U$ is at most $32$ 
and the bound is sharp.
\end{abstract}
\maketitle
\section{Introduction}
In this paper we consider three-dimensional varieties whose hyperplane
sections are Enriques surfaces. Such varieties where studied in a
number of papers started with works of G. Fano \cite{Fano-1938-e} (see
also \cite{Godeaux-1933}).

Let $U$ be a normal projective three-dimensional variety and let $A$
be a prime Cartier divisor on $U$ such that $\OOO_U(A)$ is ample. We
say that $(U,A)$ is a \textit{Fano-Enriques threefold} if $A$ is a
smooth Enriques surface and $U$ is not a generalized cone over $A$
(i.e., $U$ is not obtained by contraction of the negative section on
a $\PP^1$-bundle over $A$). Define the \textit{genus} of a
Fano-Enriques threefold by $g:=A^3/2+1$.

The main result of this paper is the following:
\begin{theorem}
The genus $g$ of a Fano-Enriques therefold $(U,A)$ is at most $17$
and the bound is sharp. Moreover, up to isomorphism there exists a
unique Fano-Enriques therefold $(U_{32},A_{32})$ of genus $17$.
\end{theorem}

Fano asserted that Fano-Enriques threefolds exists only for
$g=4,6,7,9,13$. However his arguments are unsatisfactory from the
modern point of view and, in fact, contain some gaps. In
particular, it was found latter that Fano-Enriques threefolds of
genus $g=10$ do exist. Fano's constructions were recovered by Conte
and Murre \cite{Conte-1983}, \cite{Conte-Murre-1985} but they still
assumed that singularities of $U$ are sufficiently 
general and, in fact, did
not complete the classification. New approach to the classification
problem was brought by the minimal model theory. Under additional
assumption that the singularities of $U$ are cyclic quotients,
Fano-Enriques threefolds $(U,A)$ were classified by Bayle
\cite{Bayle-1994} and Sano \cite{Sano-1995}. According to
\cite{Minagawa-1999} every Fano-Enriques threefold with only
terminal singularities admits a $\QQ$-smoothing. This means that
such a threefold is a deformation of that contained in a
Bayle-Sano's list. In particular, the genus of terminal Fano-Enriques
threefolds takes on the following values: $2\le g\le 10$ and $g=13$.
An important result was obtained by Cheltsov \cite{Cheltsov-1996a}
who proved that any Fano-Enriques threefold $(U,A)$ has only
\textup($\QQ$-Gorenstein\textup) canonical singularities. Thus to
complete the classification one has to consider the case of 
non-terminal canonical singularities.

Note that our results are stronger than that stated in the main
theorem. In fact, 
we prove many facts about Fano-Enriques threefolds,
especially, when the genus is sufficiently large.
Note also that the bound $g\le 17$ follows also from
a recent result of I. Karjemanov (in preparation).

\subsection*{Acknowledgments}
The work was carried out at Max-Planck-Institut f\"ur Mathematik,
Bonn. The author would like to thank this institute for its support
and hospitality.

\section{Preliminaries}
\label{sect-prel} 
Throughout the paper the ground field is 
supposed to be the complex number field $\CC$. 
This paper is a continuation of our previous paper
\cite{Prokhorov-2005a}, so we keep almost all the notation.
Additionally to the techniques of \cite{Prokhorov-2005a}
we use Shokurov's connectedness lemma and inversion of adjunction 
(see \cite[Ch. 17]{Utah}).
These two key results 
will be used freely without additional reference.

The following facts can be found in \cite{Conte-Murre-1985}. We only
note that there the authors assumed that $A$ is very ample and
gives an embedding of $U$ as a projectively normal variety. However
these hypothesis are not needed to prove (i)-(iii) below.
\begin{proposition}
\label{prop-predv} 
Let $(U,A)$ be a Fano-Enriques threefold of genus
$g$. Then
\begin{enumerate}
\item
$U$ has only isolated singularities,
\item
the Weil divisor $K_U$ is not Cartier,
\item
$\dim |A|=g$ and $\dim |-K_U|=g-1$.
\end{enumerate}
\end{proposition}

\begin{theorem}[\cite{Cheltsov-1996a}]
\label{th-Che} Let $(U,A)$ be a Fano-Enriques threefold. Then $U$
has only \textup($\QQ$-Gorenstein\textup) canonical singularities
and $K_X+A$ is a $2$-torsion element in the Weil divisor class
group.
\end{theorem}

\subsection{}
Let $(U,A)$ be a Fano-Enriques threefold of genus $g$. Let $\pi
\colon V\to U$ be the global log canonical cover. Then $\pi$ is a
finite \'etale in codimension two morphism such that $\pi^*(K_U+A)\sim
0$. In particular, $\pi^*K_U=K_V$ is Cartier. Hence $V$ is a Fano
threefold of degree $4g-4$ with canonical Gorenstein singularities.
Let $\tau$ be the Galous involution on $V$. Then $\tau$ acts freely
outside of a finite number of points and $U=V/\tau$.

\begin{theorem}[\cite{Jahnke-Radloff-2004},
\cite{Przhiyalkovskii-Cheltsov-Shramov-2005}] In the above notation,
if $g=g(U)\ge 12$, then the linear system $|-K_V|$ is very ample and
defines an embedding $V\subset \PP^{2g}$.
\end{theorem}

\subsection{}
\label{subs-ass-va} From now on we assume that $-K_V$ is very ample.
(However do not assume that $g\ge 12$).

\subsection{}
If the singularities of $V$ are terminal, they are isolated
cDV. Then we can apply \cite[Th. 4.2]{Minagawa-1999}.
By this result there is an one-parameter deformation $\mathfrak U\to
\mathfrak T\ni 0$ with central fibre $\mathfrak U_0=U$ and nearby
fibres having only cyclic quotient singularities. Then the total
space $\mathfrak U$ is $\QQ$-Gorenstein and $-K_{\mathfrak U_t}^3$,
$t\in \mathfrak T$ is a constant. By Bayle-Sano's classification
\cite{Bayle-1994}, \cite{Sano-1995} we have $-K_{\mathfrak U_t}^3\le
24$. Hence, $g(U)\le 13$ (and moreover $g\neq 11,\, 12$).

\subsection{}
Thus we assume that $V$ has at least one non-terminal point $P$.
We distinguish two main cases:
\begin{itemize}
\item[(I)]
$\tau P\neq P$,
\item[(II)]
the point $P$ is $\tau$-invariant.
\end{itemize}
Note that (I) is equivalent to the following:
\begin{itemize}
\item[(I${}'$)]
$U$ has a non-terminal Gorenstein point.
\end{itemize}

\subsection{}
Let $\LLL:=|-K_V|$. This linear system has two $\tau$-invariant
subsystems:
\[
\begin{array}{lll}
\LLL^+:=&\pi^*|A|,& \dim \LLL^+=g,
\\
\LLL^-:=&\pi^*|-K_U|,\quad & \dim \LLL^-=g-1.
\end{array}
\]
The following lemma is an immediate consequence of our construction.
\begin{lemma}
A general member $L^+\in \LLL^+$ is a smooth $\tau$-invariant K3
surface. Furthermore, the pair $(V,\LLL^+)$ is canonical.
\end{lemma}

\subsection{}
We will apply different kinds of the log minimal model program
(LMMP) in the category of $G$-varieties (with $G=\muu_2=\langle
\tau\rangle$). For a very brief introduction we refer to \cite[\S
2.2]{KM-1998}. In fact, there is no big difference between 
the $G$-LMMP
and standard LMMP. We only emphasize the following:
\begin{itemize}
\item
$G$-LMMP deals with $G\QQ$-factorial varieties. The latter means
that every $\tau$-invariant divisor is $\QQ$-Cartier.
\item
If we work with log divisor $K+D$, where $D$ is a boundary (resp.
linear system), this $D$ should be $G$-invariant.
\item
Instead of the Picard group $\Pic X$, Mori cone $\Mori (X)$, etc., we
should consider their $G$-invariant analogs $\Pic^G X$, $\Mori^G
(X)$, etc.
\item
Every divisorial contraction decreases the invariant Picard number
$\rho^G(X)$ by $1$ and contracts a $G$-invariant divisor.
\end{itemize}
In particular, for every threefold $V$ with canonical singularities,
one can construct a $G\QQ$-factorial terminal modification
$\phi\colon W\to V$. This is by definition a birational
$G$-equivariant contraction such that $W$ has only terminal
$G\QQ$-factorial singularities and $K_W=\phi^* K_V$. Such a
modification is not unique but every two of them are related by a
sequence of $G$-equivariant flops.

\subsection{}
\label{not-lin-syst} If $\LLL$ is a linear system on $V$ without
fixed component and $X$ is a birational model of $V$, we denote by
$\LLL_X$ the birational transform of $\LLL$ on $X$. By the abuse of
notation we often will write simply $\LLL$ instead of $\LLL_W$.

\section{An example}
\label{sect-exa} Let $x$ and $y_{i, j}$, $0\le i,\, j\le 2$ be
homogeneous coordinates in $\PP^9$. Consider the anti-canonical
embedding of $S=\PP^1\times \PP^1$ in $\PP^8=\{x=0\}\subset \PP^9$:
\[
\begin{array}{lll}
(u_0:u_1)\times (v_0:v_1) &\longmapsto& (y_{0, 0}:\dots : y_{2, 2}),
\\[8pt]
&& y_{i, j}= u_0^iu_1^{2-i}v_0^jv_1^{2-j}.
\end{array}
\]
Let $V\subset \PP^9$ be the projective cone over $S$ and let
$P=(0:\cdots:0:1)$ be its vertex.
\begin{lemma}
The variety $V$ is a Gorenstein Fano threefold with canonical
singularities. Moreover, $-K_V=2M$, where $M$ is the class of
hyperplane sections.
\end{lemma}
\begin{proof}
Since $S\subset \PP^8$ is projectively normal, $V$ is normal. Let
$\sigma \colon W\to V$ be the blowup of $P$. Then $W$ is a
$\PP^1$-bundle over $S$ and $\sigma$ contracts its negative section
$E$ to $P$. More precisely, $W\simeq \PP(\OOO_S\oplus
\OOO_{S}(-K_S))$ and the map $\sigma\colon W\to V\subset \PP^9$ is
given by the tautological linear system $|\OOO(1)|$. Hence,
$\OOO(1)\sim \sigma^*M$. Since $K_W\sim \OOO(-2)$ and this divisor
is trivial on $E$, $K_V$ is Cartier, $K_W=\sigma^*K_V$ and
$K_V=-2M$.
\end{proof}

Define the action of cyclic group $\gr=\{1,\, \tau\}$ on $\PP^8$ and
$\PP^9$ via
\[
\tau\colon x \longmapsto -x, \qquad y_{i, j} \longmapsto (-1)^{i+j}
y_{i, j}.
\]
Then both $S$ and $V$ are $\tau$-invariant and the induced action
on $S$ is
\[
(u_0:u_1)\times (v_0:v_1) \longmapsto (-u_0: u_1)\times (-v_0: v_1).
\]
The locus of $\tau$-fixed points in $\PP^9$ consists of two
projective subspaces
\begin{eqnarray*}
\PP^4_{+}&:=&\{y_{0, 0}= y_{0, 2}= y_{2, 0}= y_{2, 2}=y_{1, 1}=0\},
\\
\PP^4_{-}&:=&\{y_{0, 1}= y_{2, 1}= y_{1, 0}=y_{1, 2}= x=0\}.
\end{eqnarray*}
Clearly,
\[
\PP^4_{+}\cap V=\{P\},\qquad \PP^4_{-}\cap V=\{P_{0, 0}, P_{0, 2},
P_{2, 0}, P_{2, 2}\},
\]
where $P_{i, j}:=\{x=0, y_{k, l}=0 \mid (k, l)\neq (i, j)\}$. In
particular, the action of $\tau$ on $V$ is free in codimension two.

\begin{proposition}
\label{th-eaxa} The quotient of $V$ by $\gr=\{1,\, \tau\}$ is a
Fano-Enriques threefold $U_{32}$ of genus $17$.
\end{proposition}

\begin{proof}
Let $\LLL^+$ be the linear system that cuts out on $V$ by quadrics of
the following form
\[
q_1(y_{0, 0}, y_{0, 2}, y_{2, 0}, y_{2, 2}, y_{1, 1}) +q_2(y_{0,
1}, y_{2, 1}, y_{1, 0}, y_{1, 2}, x)=0,
\]
where $q_1$ and $q_2$ are quadratic homogeneous forms. It is easy to
see that $\LLL^+$ is base point free and each member of $\LLL^+$ is
$\tau$-invariant. In particular, a general member $L\in \LLL^+$ is
smooth and does not contain any of $P$, $P_{0, 0}$, $P_{0, 2}$,
$P_{2, 0}$, $P_{2, 2}$. Therefore the action of $\tau$ on $L$ is
free. Since $\LLL^+\subset |-K_V|$, $L$ is a K3 surface. Let
$\pi\colon V\to U=V/\tau$ be the quotient morphism and let
$A:=\pi(L)=L\tau$. 
Then $A$ is a smooth Enriques surface. Finally, we have
$L=\pi^*A$ and $2g-2=A^3=\frac12 L^3=\frac 12 (2 M)^3=32$. Hence,
$g=17$.
\end{proof}

\begin{remark}
A similar construction can be applied to anticanonically embedded
del Pezzo surface $S=S_6\subset \PP^6$ of degree $6$. We get a
new
Fano-Enriques threefold of genus $g=13$ with canonical
singularities.
\end{remark}

\section{Case (I)}
In this section we consider the case when $U$ has a non-terminal
Gorenstein point. Recall that $-K_V$ is very ample by 
our assumption (see
\ref{subs-ass-va}). We prove the following
\begin{proposition}
\label{prop-II-case-m} Under the above assumption $g\le 5$.
\end{proposition}

\begin{lemma}
\label{lemma-II-case-1} Let $(U,A)$ be a Fano-Enriques threefold.
Assume that $U$ has a non-terminal Gorenstein point $O$ \textup(in
this lemma we do not assume that $-K_V$ is very ample\textup). Then
for a general member $A^0\in |A|$ passing through $O$, the pair
$(U,A^0)$ is lc. Moreover, $O$ is the only log canonical center of
codimension $>1$.
\end{lemma}

\begin{proof}
Let $\AAA^O\subset |A|$ be the linear subsystem
consisting of divisors passing through $O$ and let
$\LLL^O:=\pi^*\AAA^O$. A general member $C$ of the ample linear
system $|\OOO_A(A)|$ is a smooth irreducible curve (see, e.g.,
\cite[Th. 4.10.2]{Cossec-Dolgachev-1989-book}). Since
$H^1(\OOO_U(A))=0$, the restriction map $H^0(\OOO_U(A))\to
H^0(\OOO_A(A))$ is surjective. Therefore, there is an element
$A^0\in \AAA^O$ such that $A^0|_A=C$. The surface $A^0$ is reduced
and has only isolated singularities. Since $A^0$ is an ample
divisor, it is connected. 

Let $\pi^{-1}(O)=\{P,\, P'\}$. Then points $P,\, P'\in V$ are
non-terminal. 
Let $L^0:=\pi^*A^0\in \LLL^O$. This
surface is also reduced, irreducible and has only isolated
singularities. Since $L^0$ is a Cartier divisor, it is normal.
Clearly, $L^0$ is Gorenstein and $K_{L^0}\sim 0$. Further, points
$P,\, P'\in L^0$ are not Du Val. Indeed, otherwise points $P,\,
P'\in V$ are isolated cDV, hence terminal. By Shokurov's
connectedness result (see \cite[12.3.1]{Utah}) the singularities of
$L^0$ are log canonical and $L^0$ has no other non-Du Val
singularities. By the inversion of adjunction 
and connectedness lemma the pair
$(V,L^0)$ is lc. Moreover, the only log canonical centers of
codimension $>1$ are points $P$ and $P'$. Now the assertion follows
by the ramification formula (see, e.g., \cite[proof of 20.3]{Utah}
or 
\cite[proof of 5.20]{KM-1998}).
\end{proof}

\begin{proof}[Proof of Proposition \xref{prop-II-case-m}]
We use notation of the proof of Lemma \xref{lemma-II-case-1}. Let
$\Delta:=\frac12 L^0+\frac12L^1$, where $L^0,\, L^1\in \LLL^O$ are
general members. Then the pair $(V,\Delta)$ is lc (see
\cite[2.17]{Utah}) and its log canonical locus coincides with
$\{P,\, P'\}$. Since $L^0$ and $L^1$ are Cartier, each crepant
exceptional divisor over $P$ or $P'$ has discrepancy $-1$ with
respect to $K_V+\Delta$.

Now we need our assumption that $-K_V$ is very ample. Let
$\HHH\subset |-K_V|$ be the linear system of hyperplane sections
passing through $P$. Since $\LLL^O\subset \HHH$, the pair $(V,\HHH)$
is lc. 
From now
on we ignore the $\tau$-action, so our constructions will not be
$\tau$-equivariant. 
Let $\phi\colon W\to V$ be a terminal $\QQ$-factorial
modification such as in \cite[Lemma 6.6]{Prokhorov-2005a}. 
We can write
\[
\phi^*(K_V+\Delta)=K_W+\Delta_W+E+E',
\]
where $E$ and $E'$ are (reduced) exceptional divisors over $P$ and
$P$, respectively. In notation of \cite[\S 6]{Prokhorov-2005a} we
also have
\[
K_W+\HHH_W+B\sim \phi^*(K_V+\HHH),
\]
where $\Supp(B)=\Supp(E)$. Since $(V,\HHH)$ is lc, $B$ is reduced.
Moreover, since $L^0,\, L^1\in \HHH$, $B=E$. Now as in \cite[\S
6]{Prokhorov-2005a} we run $(K_W+\HHH_W)$-MMP. On each step, any
fibre meets $E$. So by the connectedness lemma it does not meet $E'$. At
the end we get an extremal contraction $f\colon X\to Z$ to a
lower-dimensional variety. The log divisor $K_X+\Delta_X+E_X+E_X'$
is lc and numerically trivial. Since $E_X$ and $E_X'$ do not meet
each other and $\rho(X/Z)=1$, the only possibility is when $f$ has
one-dimensional fibres. Then according to \cite[Lemma
10.1]{Prokhorov-2005a} $Z$ is smooth rational surface and $f$ is a
$\PP^1$-bundle. In this case, divisors $E_X$ and $E_X'$ must be
disjointed sections of $f$. Then the components of $\Delta_X$ are
$f$-vertical. Let $\LLL^{P,P'}$ be the linear system of hyperplane
sections of $V$ passing through $P$ and $P'$. Let $\LLL_W^{P,P'}$
and $\LLL_X^{P,P'}$ be its birational transforms on $W$ and $X$,
respectively. Then $\dim \LLL_X^{P,P'}=\dim \LLL-2=2g-2$. Since
$\LLL^O\subset \LLL_X^{P,P'}$, we can write
$\phi^*(K_V+\LLL_W^{P,P'})=K_W+\LLL_W^{P,P'}+E+E'$. By the above
all the members of $\LLL^{P,P'}_X$ are $f$-vertical. Thus
$\LLL^{P,P'}_X=f^* \FFF$, where $\FFF$ is a linear system on $Z$
whose general member is reduced and irreducible. By the adjunction
formula we have $K_{E_X}=(K_X+E_X)|_{E_X}\sim -L^0_X|_{E_X}$. Hence
$\LLL^{P,P'}_X|_{E_X}\subset |-K_{E_X}|$. Since $f|_{E_X}\colon
E_X\to Z$ is an isomorphism, $\FFF\subset |-K_Z|$. In particular,
$-K_Z$ is nef.

If $-K_Z$ is big, then by Riemann-Roch and Kawamata-Viehweg
vanishing
\[
2g-2=\dim \LLL^{P,P'}=\dim \FFF\le\dim |-K_Z|=K_Z^2\le 9.
\]
Hence $g\le 5$. Otherwise $\dim |-K_Z|\le 1$ and as above $2g-2=
\dim \FFF\le 1$, a contradiction.
\end{proof}

\section{Case (II)}
\subsection{Construction}
Throughout this section we assume that there is a $\tau$-invariant
non-terminal point $P\in V$. Recall that $-K_V$ is very ample by our
assumption \ref{subs-ass-va}. Let $\phi\colon W\to V$ be a terminal
$\gr\QQ$-factorial modification. By definition, $\phi$ is
$\tau$-equivariant, $K_W=\phi^*K_V$, and $W$ has only terminal
(Gorenstein) $\gr\QQ$-factorial singularities. Let $\HHH$ be the
linear system of hyperplane sections of $V\subset \PP^{2g}$ passing
through $P$. Clearly, $\dim \HHH=2g-1$. Following \ref{not-lin-syst}
we denote $\HHH_W$, $\LLL_W$, and $\LLL_W^+$ birational transforms
on $W$ of $\HHH$, $\LLL$, and $\LLL^+$, respectively.

\begin{lemma}[cf. {\cite[Lemma 6.6]{Prokhorov-2005a}}]
\label{lemma-prop-MMP} The linear system $\HHH_W$ consists of
Cartier divisors and the pair $(W,\HHH_W)$ is canonical. Moreover,
one can choose $\phi$ so that $\HHH_W$ is nef.
\end{lemma}
\begin{proof}
Since the class $\HHH_W$ is $\tau$-invariant, $\HHH_W$ is
$\QQ$-Cartier. Hence $\HHH_W$ consists of Cartier divisors
(recall that $W$ has only terminal Gorenstein
singularities).
According to \cite[Lemma 6.6]{Prokhorov-2005a} there is (a
non-$\tau$-equivariant) terminal modification $\phi'\colon W'\to V$
such that $(W',\HHH_{W'})$ is canonical. Varieties $W$ and $W'$ are
isomorphic in codimension one over $V$. Hence the corresponding map
$(W,\HHH_{W})\dashrightarrow(W',\HHH_{W'})$ is crepant and
$(W,\HHH_W)$ is also canonical. The last statement is proved similar to
the corresponding statement in \cite[Lemma 6.6]{Prokhorov-2005a}.
\end{proof}

Now starting with $(W,\HHH_W)$ we run $\gr$-LMMP with respect to
$K+\HHH$. Similar to \cite[Lemma 3.4]{Prokhorov-2005a} we see that
the properties in Lemma \ref{lemma-prop-MMP} are preserved. 
At the end we get a $\gr$-extremal $K+\HHH$-negative
contraction $f\colon X\to Z$ to a lower-dimensional variety. By the
above and inversion of adjunction we have

\begin{corollary}
\label{corollary-H-X} The pair $(X,\HHH)$ is canonical. Hence a
general member $H\in \HHH$ is a normal surface with at worst Du Val
singularities. Moreover, $H$ does not pass through non-Gorenstein
points of $X$.
\end{corollary}

\begin{lemma}
\label{lemma-action-mu-X} 
For a general member $L\in \LLL^+_X$, the
quotient $L/\tau$ is an Enriques surface with at worst Du Val
singularities. In particular, the action of $\tau$ on $L$ is free.
\end{lemma}
\begin{proof}
Since all birational transformations $V\dashrightarrow X$ are
$K+\LLL^+$-crepant, the pair $(X,L)$ is canonical and $K_X+L\sim 0$.
Hence $L$ is a K3 surface with at worst Du Val singularities. The
surface $L/\tau$ has only quotient singularities and is birationally
equivalent to an Enriques surface. The action of $\tau$ extends to
the minimal resolution $\tilde L$ of $L$. We have the following
diagram
\[
\xymatrix{ \tilde L\ar[r]^{\tilde \varsigma}\ar[d] &\tilde
L/\tau\ar[d]^{\upsilon}
\\
L\ar[r]^{\varsigma}&L/\tau }
\]
By the ramification formula $K_{\tilde L}=\tilde
\varsigma^*K_{\tilde L/\tau} +\tilde \Xi$, where $\tilde \Xi$ is the
branch divisor. Since $K_{\tilde L}=0$ and $\tilde L/\tau$ is
birationally equivalent to an Enriques surface, we have $\tilde
\Xi=0$ and $K_{\tilde L/\tau}\equiv 0$. It follows that the action
of $\tau$ on $\tilde L$ and $L$ is free in codimension one and
$\upsilon$ is crepant: $K_{\tilde L/\tau}=\upsilon^* K_{L/\tau}$.
Hence the singularities of $L/\tau$ are Du Val and $K_{L/\tau}\equiv
0$. This proves the first statement. Further, for the topological
Euler number we have $24=\chi_{\operatorname{top}}(\tilde L)=
2\chi_{\operatorname{top}}(\tilde L/\tau)-s$, where $s$ is the
number of $\tau$-fixed points. This gives us
$\chi_{\operatorname{top}}(\tilde L/\tau)= 12+s/2$. Since $\tilde
L/\tau$ is an Enriques surface with singularities of type $A_1$, we
have $s=0$.
\end{proof}

Now we consider cases according to the dimension of $Z$. The main
theorem is a consequence of Propositions \ref{dimZ=0}, \ref{dimZ=1},
and \ref{prop-surf-main} below.

\subsection{Case: $\dim Z=0$}
Then $\rho^{\langle\tau\rangle}(X)=1$, so $-K_X\equiv rH$, where
$H\in \HHH_X$ and $r\in \QQ$. Since $-(K_X+H)$ is ample, $r>1$. By
\cite[\S 7]{Prokhorov-2005a} we have $\dim |-K_X|\le 33$. Therefore,
$g\le 16$ in this case. Analysing the properties of action of $\tau$ we
can get more precise result:

\begin{proposition}
\label{dimZ=0} If $Z$ is a point, then $g\le 5$.
\end{proposition}
\begin{proof}
By the adjunction formula $-K_H\equiv (r-1)H|_H$ is ample. Hence $H$
is a del Pezzo surface with at worst Du Val singularities. So,
$K_H^2\le 9$. By Kawamata-Viehweg vanishing $h^1(\OOO_X)=0$.
Therefore, by  Riemann-Roch 
on $H$ we have
\[
2g\le h^0(\OOO_X(H))=h^0(\OOO_H(H))+h^0(\OOO_X)=
\frac{r}{2(r-1)^2}K_H^2+2.
\]
If $X$ is Gorenstein, then $r\ge 2$ and $g\le 5$.

Assume that $X$ is not Gorenstein. Then according to
\cite{Sano-1996} $X$ is either a weighted projective space
$\PP(1,1,1,2)$ or isomorphic to one of the following weighted
hypersurfaces:
\begin{itemize}
\item
$X_6\subset\PP(1,1,2,3,I)$, $I=2,3,4,5,6$;
\item
$X_4\subset\PP(1,1,1,2,I)$, $I=2,3$;
\item
$X_3\subset\PP(1,1,1,1,2)$.
\end{itemize}
If $X\simeq \PP(1,1,1,2)$, then $X$ has a unique point
$P:=(0:0:0:1)$ of index $2$. This point is contained in the base
locus of $\LLL^+\subset |-K_X|$. On the other hand, $P$ must be
$\tau$-invariant. This contradicts Lemma \ref{lemma-action-mu-X}.
Similar arguments show that cases $X_6\subset\PP(1,1,2,3,4)$,
$X_6\subset\PP(1,1,2,3,5)$, $X_4\subset\PP(1,1,1,2,3)$, and
$X_3\subset\PP(1,1,1,1,2)$ are also impossible. In case
$X_6\subset\PP(1,1,2,3,2)$ the curve $\Upsilon$ given by
$x_1=x_2=x_4=0$ is $\tau$-invariant. The intersection $X_6\cap
\Upsilon$ is given on $\Upsilon =\PP(2,2)\simeq\PP^1$ by a cubic
polynomial. Hence there is a $\tau$-invariant point $P\in X_6\cap
\Upsilon$. This point is of index $2$, so $P\in \Bs |-K_X|$. As
above we get a contradiction. Similar arguments work in the case
$X_6\subset\PP(1,1,2,3,6)$. It remains to consider two cases:
\begin{itemize}
\item
$X\simeq X_6\subset\PP(1,1,2,3,3)$, and
\item
$X\simeq X_4\subset\PP(1,1,1,2,2)$.
\end{itemize}
In the first case, we have $\OOO_X(H)\simeq \OOO_X(3)$. Then $2g\le
h^0(\OOO_X(H))=4+4+2=10$. The second case is treated similarly.
\end{proof}

\subsection{Case: $\dim Z=1$}
Then $Z$ is a smooth rational curve.
\begin{proposition}
\label{dimZ=1} 
If $Z$ is a curve, then $g\le 16$.
\end{proposition}
Since $\rho^{\langle\tau\rangle}(X/Z)=1$, we have 
$-K_X\equiv r H+f^*\Xi$, where $r\in \QQ$, $r>1$,
and $\Xi$ is a $\QQ$-divisor on $Z$.
By adjunction formula the generic fibre $X_\eta$ 
is isomorphic either $\PP^2$ or $\PP^1\times \PP^1$.
We need some information about the structure of
singular fibres.

\begin{lemma}
\label{lemma-e-dp-f} 
Let $f\colon X\to Z$ be a contraction from a
terminal threefold to a curve. Assume there exists an $f$-ample
Cartier divisor $H$ such that $-(K_X+H)$ is $f$-ample and $(X,H)$ is
canonical. Then every fibre $f^*o$, $o\in Z$ is reduced. Moreover,
$f^*o$ has at most two components
and the pair $(X,f^*o)$ is dlt.
If $f^*o$ is irreducible, then it is isomorphic to
one of the following surfaces:
\begin{equation}
\label{eq-delta-genus}
\PP^2,\quad \PP(1,1,4),\quad \PP^1\times \PP^1,
\quad \PP(1,1,2).
\end{equation}
\end{lemma}
Note that proofs of corresponding statements in
\cite{Prokhorov-2005a} (Corollary 8.6 and Lemma 9.2) have some small
gaps. The above lemma is a generalization and correction of these
statements.
\begin{proof}
As above the generic fibre $X_\eta$ 
is isomorphic either $\PP^2$ or $\PP^1\times \PP^1$.
Fix a point $o\in Z$ and
regard $Z$ and $X$ as small neighborhoods of $o$ and the
fibre $f^{-1}(o)$, respectively. 
By the inversion of adjunction the surface $H$ has only Du Val singularities 
and does
not contain any non-Gorenstein points of $X$. 
The restriction $\varphi\colon H\to Z$ is a rational curve fibration
such that $-K_H$ is $\varphi$-ample.
Since $K_H$ is Cartier, every fibre of $\varphi$ is 
isomorphic to a plane conic.
In particular, so is $\varphi^*o=f^*o\cap H$. 
This immediately implies that every ample divisor
on $H$ is very ample (over $Z$).
Further, by Kawamata-Viehweg vanishing
$H^1(\OOO_X)=0$. Hence the restriction map $H^0(\OOO_X(H))\to
H^0(\OOO_H(H))$ is surjective and 
$|H|$ is $f$-base point free.

If the fibre $f^*o$ is not reduced, $\varphi^*o$ is a
double line. So $f^*o=2S$, where $S$ is a prime Weil $\QQ$-Cartier
divisor and $S\cap H$ is a reduced irreducible rational curve. Since
$X$ is Gorenstein (and terminal) near $H$, 
$S\cap H$ is Cartier on $H$. Thus $H$ is
smooth near $S\cap H\simeq \PP^1$. 
But then the rational curve fibration $\varphi\colon H\to
Z$ cannot have multiple fibres, a contradiction.
Therefore, $f^*o$ is reduced.

If the fibre $S:=f^*o$ is irreducible, 
all the arguments \cite[Lemmas 8.1, 8.5]{Prokhorov-2005a} work.
Indeed, $\chi(\OOO_S(H))=\chi(\OOO_{X_\eta}(H))$
and $h^2(\OOO_S(H))=0$. Hence, 
$h^0(\OOO_S(H))\ge h^0(\OOO_{X_\eta}(H))$.
For the $\Delta$-genus, we have 
\[
\Delta(S,\OOO_S(H)):=
\dim S+\deg \OOO_S(H)- h^0(\OOO_S(H))
\le \Delta(X_\eta,\OOO_{X_\eta}(H))=0.
\]
By Fujita's classification of polarized varieties of 
$\Delta$-genus zero we obtain the cases in \eqref{eq-delta-genus}.
By the inversion of adjunction $(X,S)$ is plt in these cases.

From now on we assume that $f^*o$ has at least 
two components. Since $|H|$ is base point free, 
replacing $H$ with a general member of $|H|$
we also may assume that $\varphi^*o$ is reduced
(and reducible). 
Thus $\varphi^*o$ is isomorphic to a pair of lines 
meeting in a point $P$: $\varphi^*o=l+l'$, $l_i\simeq\PP^1$, 
$l\cap l'=\{P\}$.
In particular, this implies that $f^*o$ 
has exactly two components:
$f^*o=S+S'$.
Note that $S\cap S'$ is of pure 
dimension one (see for example \cite[Lemmas 8.7]{Prokhorov-2005a}).
We claim that $(H,\varphi^*o)$ is lc.
Indeed, by adjunction
\[
(K_H+l+l')|_{l}=K_{l}+\Diff_{l}(l').
\]
Hence $\deg \Diff_{l}(l')=2+K_H\cdot l<2$.
Clearly, $\Supp \Diff_{l}(l')=P$.
Since $K_H+l+l'$ is Cartier, $\Diff_{l}(l')$
is an integral divisor. Thus, $\Diff_{l}(l')=P$
and the pair $(l,\Diff_{l}(l'))$ is lc.
By the inversion of adjunction so is $(H,\varphi^*o)$.

Again 
by the inversion of adjunction the pair
$(X,H+f^*o)$ is lc near $H$. 
By Shokurov's connectedness lemma $(X,H+f^*o)$ is lc 
everywhere.
Let $\Gamma$ be the locus of log canonical singularities 
of the pair $(X,f^*o)$. Then $\Gamma\cap H=\{P\}$.
Since $|H|$ is base point free, 
the one-dimensional component $\Gamma_0\subset \Gamma$ 
is an irreducible curve.
Assume that $(X,f^*o)$ has a
zero-dimensional center $P$ of log canonical singularities. 
By the above
$P\notin H$. Take a general hyperplane section $M$ passing through
$P$. Then, for any $\epsilon>0$, the pair $(X,f^*o+H+\epsilon M)$ is
not lc at $P$ and lc near $H$. Again applying connectedness
lemma to $(X,(1-\delta)f^*o+H+\epsilon M)$ for $0<\delta \ll
\epsilon\ll 1$ we derive a contradiction. Thus $(X,f^*o)$ 
has no any
zero-dimensional centers of log canonical singularities
and its log canonical locus is an irreducible curve
$\Gamma=S\cap S'$. 
Hence $(X,f^*o)$ is dlt (see \cite[Prop. 2.40]{KM-1998}).
\end{proof}

Now by \cite[16.15]{Utah} we have the following two corollaries.
\begin{corollary}
In notation of Lemma \xref{lemma-e-dp-f} each component of $f^*o$ is
normal.
\end{corollary}

\begin{corollary}
\label{cor-tw-comp-sing} Notation as in Lemma \xref{lemma-e-dp-f}.
Assume that the fibre $f^*o$ is reducible and let $S,\, S'$
be its irreducible components.
Let $P\in S\cap S'\subset X$ be a point of index $m>1$.
Then, near $P$, each component $S,\, S'\subset
f^*o$ is analytically isomorphic to $\CC^2/\muu_m(1,q_i)$, $\gcd
(q_i,m)=1$ and the intersection curve $S\cap S'$ is smooth.
\end{corollary}

\begin{lemma}
\label{lemma-fiber-surf} Let $S$ be a normal surface. Assume that
$-K_S\equiv \alpha h +C$, where $\alpha>1$, $h$ is an ample Cartier
divisor, and $C$ is an effective Weil divisor. Then $S$ has at most
one singular point. If $\alpha\ge 2$, then $S\simeq \PP^2$.
\end{lemma}
\begin{proof}
Let $\mu\colon \tilde S\to S$ be the minimal resolution, let $\tilde
C$ be the proper transform of $C$, and let $h^*=\mu^*h$. We can
write
\[
K_{\tilde S}\equiv \mu^*K_{S}+\Delta, \quad \mu^*C=\tilde C+\Theta,
\]
where $\Delta$ and $\Theta$ are $\mu$-exceptional effective
$\QQ$-divisors. Therefore,
\[
K_{\tilde S}+\alpha h^*\equiv -(\tilde C+\Delta+\Theta).
\]
Let $R$ be a $K_{\tilde S}+\alpha h^*$-negative extremal ray on
$\tilde S$, let $\varphi\colon \tilde S\to Q$ be the
corresponding contraction, and let $\ell$ be a curve contracted by
$\varphi$. Then $h^*\cdot \ell\ge 1$, $\Delta\cdot \ell\ge 0$, and
$\Theta\cdot \ell\ge 0$. Hence $K_{\tilde S}\cdot \ell\le
-\alpha<-1$. There are two possibilities:
\begin{enumerate}
\item
$Q$ is a point. Then $S=\tilde S\simeq \PP^2$.
\item
$Q$ is a curve and $\varphi$ is a $\PP^1$-bundle. Then $\rho(\tilde
S)=2$. 
\end{enumerate}
In the first case we are done. Consider the second case.
Since $\rho(\tilde
S)=2$, $\mu$ either is an identity map or contracts an
irreducible curve. Hence $S$ has at most one singular point.
Assume
$\alpha\ge 2$. Since $K_{\tilde S}\cdot \ell=-2$, $h^*\cdot \ell=
1$. Then $\alpha=2$ and $(\tilde C+\Delta+\Theta)\cdot \ell=0$. This
is possible only if $\Delta=\Theta=0$ and $\tilde C$ is contained in
fibres of $\varphi$. If $\mu$ is not an identity map, then $\tilde
C$ meets the $\mu$-exceptional curve. Then $\Theta\neq 0$, a
contradiction. Hence $S$ is smooth and $-K_{S}\equiv 2h+ C$ is
ample. Thus there is a $K$-negative extremal ray $R'\neq R$. For
corresponding extremal curve $\ell'$ we have $C\cdot \ell'>0$, so
$-K_S\cdot \ell'\ge 3$, a contradiction.
\end{proof}

\begin{corollary}
\label{cor-1-n-point}
Assume that $f^*o$ is reducible. Then
\begin{enumerate}
\item
If $X_\eta\simeq\PP^2$, then $X$ has a unique non-Gorenstein point
on $f^*o$.
\item
If $X_\eta\simeq\PP^1\times\PP^1$, then $X$ is Gorenstein near
$f^*o$.
\end{enumerate}
\end{corollary}
\begin{proof}
Let $f^*o=S+ S'$. Near $f^*o$ we have $-K_X\equiv \alpha H$,
where $\alpha=3/2$ (resp., $\alpha=2$) in the case
$X_\eta\simeq\PP^2$ (resp., in the case
$X_\eta\simeq\PP^1\times\PP^1$). By the adjunction formula we have
\[
-K_{S}\equiv - \left(K_X+S\right)|_{S}\equiv \left(\alpha H
+S'\right)|_{S}\equiv \alpha h +C,
\]
where $h:=H|_{S}$. If $X_\eta\simeq\PP^1\times\PP^1$, then by Lemma
\ref{lemma-fiber-surf} we have $S\simeq\PP^2$ and by symmetry
$S'\simeq\PP^2$. By Corollary \ref{cor-tw-comp-sing} $X$ is
Gorenstein near $f^*o$ in this case. Consider the case
$X_\eta\simeq\PP^2$. If $X$ is Gorenstein near $f^*o$, then
$K_X^2\cdot S=1/2 K_X^2\cdot f^*o=9/2$ must be an integer, a
contradiction. Then the assertion follows by Lemma
\ref{lemma-fiber-surf} and Corollary \ref{cor-tw-comp-sing}.
\end{proof}

\begin{proof}[Proof of Proposition \xref{dimZ=1}]
Assume $g\ge 17$.
Consider the case when $X_\eta$
is isomorphic to $\PP^2$. 
Let $o\in Z$ be a $\tau$-fixed point and let $S:=f^*o$ be the fibre.
If $S\simeq \PP^2$, then the locus of $\tau$-fixed
points on $S$ consists of a line $\Gamma$ and an isolated point. But
then each member of $\LLL^+$ meets $\Gamma$.
This contradicts Lemma \ref{lemma-action-mu-X}.
If $S\not \simeq \PP^2$, then by Lemma \ref{lemma-e-dp-f} 
and Corollary \ref{cor-1-n-point} 
$X$ has only one 
non-Gorenstein point $P\in f^*o$.
This point must be $\tau$-invariant.
On the other hand,
$P\in \Bs |-K_X|\subset 
\Bs \LLL^+$. Again we have a contradiction
by Lemma \ref{lemma-action-mu-X}.

Now we consider the case when $X_\eta$
is isomorphic to $\PP^1\times \PP^1$.
By Lemma \ref{lemma-e-dp-f} 
and Corollary \ref{cor-1-n-point} $X$ is Gorenstein
and each fibre of $f$ is a reduced quadric in $\PP^3$.
In this case the linear system $|H|$ determines a
$\tau$-equivariant
embedding into $\PP:=\PP(\EEE)$, where $\EEE$ is a rank $4$
vector bundle over $Z\simeq \PP^1$. We may assume that 
$\EEE=\oplus \OOO(d_i)$, where $d_1\ge d_2\ge d_3\ge d_4=0$.
As in \cite[\S 9]{Prokhorov-2005a} we introduce 
the following notation.
Put $d=\sum d_i$.
Let $M$ and $F$ be classes of the tautological divisor
and a fibre of the projection $\PP\to Z$.
Then $X\sim 2M+rF$ for some $r\in \ZZ$.
Let $G$ and $Q$ be restrictions on $X$ of $M$ and $F$,
respectively. 
Clearly, $\Pici X=
\ZZ\cdot G\oplus \ZZ\cdot Q\simeq \Pic \PP$.
Hence, 
\[
H\sim G+\alpha Q,\quad B\sim G-(d+r+\alpha-2)Q,
\quad \alpha\in \ZZ.
\]
Since $R^if_*\OOO_X(H)=R^if_*\OOO_X(B)=0$
for $i>0$,
\[
H^0(\OOO_X(H))\simeq H^0(\EEE(\alpha)),
\quad
H^0(\OOO_X(B))\simeq H^0(\EEE(2-d-r-\alpha)).
\]
Since $B$ is effective, $d+r+\alpha-2\le d_1$.
Further,
\[
L\cdot B\cdot H=(-K_X)\cdot B\cdot H=2(6-d-2r)\ge 0
\]
(because $H$ is nef, see \cite[Lemma 9.8]{Prokhorov-2005a}).
Hence, $d+2r\le 6$. Further,
\begin{multline}
\label{eq-quadrics}
2g-1\le h^0(\OOO_X(H))=d+4+4\alpha\le 
d+4
+4(2+d_1-d-r)
\\
=12+4d_1-3d-4r\le 12+d_1-4r.
\end{multline}
We claim that $r<0$. Indeed, assume $r\ge 0$.
Then $d_1\le d\le 6$, $2g-1\le 18$, and $g\le 9$,
a contradiction.

Let $C\subset \PP(\EEE)$ be the subscrol 
$\PP(\oplus_{d_i=0}\OOO(d_i))$.
The linear system $|M|$ is base point free and contracts $C$.
Since $X\sim 2M+rF$, where $r<0$, $C$ is contained in $X$.
Clearly, $\dim C=1$.
In particular, $d_2\ge 1$ and $d_1\le d-2$. 
Thus \eqref{eq-quadrics} can be rewritten as follows
\begin{equation}
\label{eq-quadrics-1}
2g-1\le d+4+4\alpha\le 
12+4d_1-3d-4r\le 4+d-4r.
\end{equation}

If $K_X\cdot C>0$, then
$C$ is contained in the base locus of $\LLL^+$.
On the other hand, $\tau$ has two fixed points on $C\simeq\PP^1$.
This contradicts Lemma \ref{lemma-action-mu-X}.
Hence $K_X\cdot C<0$ and $-K_X$ is nef
and big. From the main result of 
\cite{Prokhorov-2005a} we get 
\begin{multline*}
72\ge -K_X^3=(2G+(2-d-r)Q)^3=
(2M+rF)\cdot (2M+(2-d-r)F)^3=
\\
=16d+24(2-d-r)+8r=8(6-d-2r).
\end{multline*}
Hence, $-3\le d+2r\le 6$. Since $G\cdot C=0$, we have
$0\le -K_X\cdot C=2-d-r$. This gives us 
$d+r\le 2$, $d\le 7$, and $r\ge -5$.
By \eqref{eq-quadrics-1} we get $g\le 16$.
\end{proof}

\subsection{}
\label{not-surf-main} Consider the case when $\dim Z=2$. Then
$\HHH_X$ is a linear system of sections of $f$. Since the invariant
Picard number of $X$ over $Z$ is equal to $1$, $f$ is
equi-dimensional. In this situation, the proof of \cite[Lemma
10.1]{Prokhorov-2005a} works and we get that both $X$ and $Z$ are
smooth and $f$ is a $\PP^1$-bundle. Thus $X=\PP(\EEE)$, where $\EEE$
is a rank $2$ vector bundle over $Z$. By construction the pair
$(X,\LLL\subset |-K_X|)$ is canonical. There is a decomposition
$-K_X\sim H+B$, where $H$ and $B$ are effective $f$-ample divisors
(sections of $f$), the divisor $H$ is nef, and the image of the map
$\Phi_{|H|}$ given by the linear system $|H|$ is three-dimensional.
We have
\[
2g-1=\dim \HHH \le \dim |H|.
\]

\begin{proposition}
\label{prop-surf-main} In the above notation we have $2g-1\le \dim
|H|\le 33$. Moreover, if $\dim |H|= 33$, then $g(U)=17$ and $U\simeq
U_{32}$ \textup(see \S \xref{sect-exa}\textup).
\end{proposition}

\subsection{}
\label{lemma-action-mu-X-2}
The involution $\tau$ acts on $Z$ effectively. Let $L\in \LLL_X^+$
be a general member and let $f_L\colon L\to \bar L\to Z$ be the
Stein factorization. Here $\bar L\to L$ is a finite of degree two
morphism with branch divisor $\Theta\in |-2K_Z|$. The involution
$\tau$ naturally acts on $L$, $\bar L$ and $\Theta$. 
Moreover, by Lemma \ref{lemma-action-mu-X} $\tau$ has
no fixed points on $\Theta$.

\subsection{}
Now we run $\gr$-MMP on $Z$: $Z=Z_1\to Z_2\to \cdots \to Z_N=Z'$.
Each step is a contraction of an $\tau$-invariant set of disjointed
$(-1)$-curves. At the end we get one of the following cases
\cite{Iskovskikh-1979s-e}:
\begin{enumerate}
\item
$\Pici Z'\simeq \ZZ$ and $Z'$ is a del Pezzo surface,
\item
$\Pici Z'\simeq \ZZ\oplus \ZZ$ and there exists a contraction
$Z'\to \PP^1$. Each fibre is isomorphic to a reduced plane conic.
\end{enumerate}
According to \cite[Lemmas 5.4, 10.4]{Prokhorov-2005a} be can
construct a sequence of $\tau$-equivariant birational
$K+\LLL$-crepant transformations
\[
\xymatrix{ X=X_1\ar[d]^{f}\ar@{-->}[r]&X_2\ar[d]^{f_2}\ar@{-->}[r]
&\cdots\ar@{-->}[r]&X_N\ar[d]^{f_N}
\\
Z=Z_1\ar[r]&Z_2\ar[r]&\cdots\ar[r]&Z_N }
\]
where each square is one of transformations (iii)-(v)
of \cite[Lemma 5.4]{Prokhorov-2005a} over either a $\tau$-invariant
$(-1)$-curve or $\tau$-invariant pair of disjointed $(-1)$-curves. 
These transformations preserve all the properties of $X$, $\LLL$,
$\LLL^+$, and $\HHH$ (except for the canonical property of
$(X,\HHH)$). Moreover, by construction all the transformations are
$\tau$-equivariant and the action of $\tau$ on a general member of
$\LLL^+$ is free (see Lemma \ref{lemma-action-mu-X}). 
Replacing $X/Z$ with $X_N/Z_N$ we may assume that
$Z$ satisfies one of the conditions (i) or (ii). Let $\Omega\subset
Z$ be the locus of $\tau$-fixed points.

Consider case $\Pici Z\simeq \ZZ$. 
Then the curve $\Theta$ is ample. Since
$\Theta\cap \Omega=\emptyset$, $\Omega$ is finite. 
If $Z$ contains no $(-1)$-curves, 
then $Z\simeq \PP^2$ or $\PP^1\times \PP^1$.
It is easy to see that in both cases 
either the set of $\tau$-fixed points
contains a curve or $\Pici Z\not \simeq \ZZ$. 
Therefore, $Z$ contains a $(-1)$-curve $E$
and $K_Z$ is not divisible, i.e., $K_Z$ generates $\Pici Z$.
Since $\Pici Z\simeq \ZZ$, $E':=\tau(E)\neq E$ and 
$E'\cap E\neq \emptyset$. Moreover, $E+E'\sim -aK_Z$
for some positive integer $a$. Then
\[
a K_Z^2=-K_Z\cdot (E+E')=2.
\]
The quotient $Y:=Z/\tau$ is a del Pezzo surface with Du Val
singularities of type $A_1$. 
In particular, $K_Y$ is Cartier and $K_Y^2$ is an integer.
Hence $K_Z^2=2K_Y^2$ is even.
Thus we have the only possibility $K_Z^2=2$, $a=1$,
$E+E'\sim -K_Z$, and $E\cdot E'=2$.
Thus for every $(-1)$-curve $E$ its image $\tau(E)$ 
is uniquely defined by $\tau(E)\sim -K_Z-E$.
Recall that the Geiser involution $\tau_g$ on the del Pezzo
surface $Z$ degree $2$ is the Galois involution of 
the anticanonical double cover $Z\to \PP^2$.
For this involution we also have $\tau_g(E)\sim -K_Z-E$.
Hence the actions of $\tau$ and $\tau_g$ on $\Pic Z$
coincide.
This obviously implies $\tau=\tau_g$.
On the other hand, the Geiser involution has a curve of fixed point,
a contradiction.

Now we consider the case when $\Pici Z\simeq \ZZ\oplus \ZZ$. Let
$\psi\colon Z\to \PP^1$ be the $\tau$-equivariant conic bundle. 
Assume that the morphism 
$\psi$ is not smooth. Let $r$ be the number
of degenerate fibres. By \cite{Iskovskikh-1979s-e}
we have $K_Z^2\le 5$. 
Therefore, $r=8-K_Z^2\ge 3$. 
The involution $\tau$ interchange components of 
singular fibres. Hence every singular fibre is 
$\tau$-invariant and there are at least three $\tau$-fixed
points on $\PP^1$. This is possible only if 
the action of $\tau$ on $\PP^1$ is trivial.
In this case, the set $\Omega$ 
contains no components of fibres and meets each smooth 
fibre at two points and each singular fibre at 
one point. Hence $\Omega$ is a connected curve
such that the restriction $\psi|_\Omega$ is of degree $2$
ramified exactly over the intersection of
$\Omega$ with singular fibres.
Note also that $\Omega$ is the branch divisor 
of the quotient map $Z\to Z/\tau$. So, $\Omega$ is smooth
and irreducible.
By the Hurwitz formula,
$2p_a(\Omega)-2=r-4$. 
By \cite{Iskovskikh-1979s-e} 
the group $\Pici Z/\PP^1=1$ is generated by 
$K_Z$ and the class of the fibres.
Hence we can write
$\Omega \sim -K_Z+m\psi^*P$, where $P\in \PP^1$ is a point and $m\in
\ZZ$. Since $\Omega\cap \Theta=\emptyset$, we have $K_Z\cdot \Omega=0$.
Then 
\[
0=-K_Z\cdot\Omega=-K_Z\cdot (-K_Z+m\psi^*P)=K_Z^2+2m,
\]
\[
K_Z^2+4m=(-K_Z+m\psi^*P)^2= \Omega^2=(K_Z+\Omega)\cdot \Omega =
2p_a(\Omega)-2.
\]
Using this by Noether formula we obtain
\[
8-r=K_Z^2=2-2p_a(\Omega)=4-r,
\]
a contradiction.

Finally, we assume that $\psi\colon Z\to \PP^1$ is a smooth
morphism. Then $Z\simeq\FF_e$,
$e\ge 0$, $e\neq 1$. 
If the action of $\tau$ on $\PP^1$ is trivial,
then as above 
$\Omega$ is a curve meeting each fibre at two points.
Therefore, $\Omega$ is a disjointed union of two sections.
But in this case at least one of these sections
meets $\Theta\sim -2K_Z$. Thus the action of $\tau$ on 
$\PP^1$ is nontrivial. Since $\Omega\cap \Theta=\emptyset$,
$\Omega$ has no $\psi$-vertical components.
Therfore, $\Omega$ is a finite set.
We claim that $e=0$. Indeed, assume
that $e\ge 2$. Let $\Sigma$ be the negative section of the ruling.
Then $\Theta\cdot \Sigma=-2K_Z\cdot \Sigma=2(2-e)$. If $e>2$, then
$\Sigma$ is a component of $\Theta$. On the other hand, $\Sigma$ is
$\tau$-invariant. But then $\tau$ must have two fixed points on
$\Sigma\subset \Theta$.
This contradicts \ref{lemma-action-mu-X-2}.
Therefore $e=2$. As above,
$\Sigma$ is not a component of $\Theta$. Hence, $\Theta\cap
\Sigma=\emptyset$. Let $F:=f^{-1}(\Sigma)$. Clearly, $F\simeq \FF_m$
for some $m\ge 0$. Let $\Sigma'\subset F$ 
be a minimal section of the ruling.
Then
\[
m-2=K_F\cdot \Sigma'= K_X\cdot \Sigma' +F\cdot\Sigma'= K_X\cdot
\Sigma'-2.
\]
Thus $K_X\cdot \Sigma'=m$. If $m>0$, then $K_X\cdot \Sigma'<0$ and
$\Sigma'$ is contained in the base locus of $\LLL_X\subset |-K_X|$.
On the other hand, $\Sigma'$ is $\tau$-invariant and $\tau$ has two
fixed points on $\Sigma'$.
This contradicts Lemma \ref{lemma-action-mu-X}.
Therefore, $m=0$,
$F\simeq \PP^1\times \PP^1$, and $\Sigma'\cdot K_X=0$. This implies
that $F$ does not intersect a general member of $\LLL_X$. In
particular, $F\cap \Bs \LLL_X=\emptyset$. This means that the map
$\Psi\colon X\dashrightarrow V$ given by $\LLL_X$ is a morphism near
$F$. By the above $\Psi$ contracts $F$ to a curve. Since $F\cdot
\Sigma'=\Sigma^2=-2$, $V$ is singular along $\Psi(F)$.
This contradicts Proposition \ref{prop-predv}, (i).
Thus $e=0$ and $Z\simeq \PP^1\times \PP^1$.

For simplicity we put $c_i:=c_i(\EEE)$, $i=1,2$. Then
\[
-K_X= 2 H +f^*(-K_Z-c_1),
\]
\[
H^2=H\cdot f^*c_1-f^*c_2,\qquad H^3=c_1^2-c_2.
\]

Let $\Sigma$ and $l$ be generators of $\PP^1\times \PP^1$. We may
assume that $c_1=a\Sigma+bl$ and $c_2=c$, where $a,b,c$ are
integers. By Riemann-Roch we have:
\begin{equation*}
\chi(\EEE)=\frac12(c_1^2-2c_2-K_Z\cdot c_1)+2= (a+1)(b+1)-c+1.\eqno
(*)
\end{equation*}
Since $\EEE$ is nef, $c_1=a\Sigma+bl$ is a nef class, so $a,\, b\ge
0$. Up to permutation we may assume that $a\le b$.

\begin{lemma}[{\cite[Lemma 10.6]{Prokhorov-2005a}}]
\label{lemma-p1-p2-nachalo} Let $\Gamma\subset Z$ be a smooth
rational curve such that $\dim |\Gamma|>0$ and let
$\EEE|_\Gamma\simeq \OOO_{\PP^1}(d_1)+\OOO_{\PP^1}(d_2)$. Then
$|d_1-d_2|\le 2+\Gamma^2$.
\end{lemma}

\begin{lemma}
\label{lemma-decomp} If $\EEE$ is decomposable, then $g\le 17$.
Moreover, equality holds only if $\EEE\simeq \OOO(4,2)\oplus
\OOO(2,4)$ and in this case $U$ is isomorphic to 
the variety $U_{32}$ from \S
\xref{sect-exa}.
\end{lemma}
\begin{proof}
Let $\EEE=\EEE_0\oplus\EEE_1$, where $\EEE_0=\OOO(\alpha,\beta)$,
$\EEE_1=\OOO(\gamma,\delta)$. We have
\[
\begin{array}{l}
a=\alpha+\gamma\ge 0,\qquad
b=\beta+\delta\ge 0, 
\\[10pt]
c_2=\alpha\delta+\beta\gamma,\quad
|\alpha-\gamma|\le 2,\quad |\beta-\delta|\le 2, 
\end{array}
\]
(by Lemma
\ref{lemma-p1-p2-nachalo}). Put $r:=\alpha-\gamma$ and
$s:=\beta-\delta$. Then $c_2=\frac12(ab-rs)$ and
\begin{equation*}
\label{eq-dim-Fn} 2g\le \chi(\EEE)= \frac12ab+a+b+\frac12rs+2=
\frac12(a+2)(b+2)+\frac12rs.
\end{equation*}
Since $-2\le r, s\le 2$, we have $g\le 17$. The equality holds only
if $a=b=6$ and $r=s=\pm 2$. Thus we may assume that $\EEE\simeq
\OOO(4,2)\oplus \OOO(2,4)$. Then $-K_X=2(H-2f^*(\Sigma+l))$, $\dim
\HHH=33$, and $\dim \LLL=34=\dim |-K_X|$. In this case, $V$ is the
anti-canonical image of $X=\PP(\EEE)=\PP(\OOO\oplus\OOO(2,2))$.
Clearly, this image coincides with the variety $V$ from \S
\ref{sect-exa}. There is only one choice for the action of $\tau$.
Hence we get the Fano-Enriques threefold $U_{32}$ from Proposition
\ref{th-eaxa}.
\end{proof}

From now on we assume that $\EEE$ is indecomposable. 
Our arguments below are very similar to that in 
\cite[\S 10]{Prokhorov-2005a}.
Let $B_0$ be
the $f$-horizontal component of $B$. Since $L\cap B$ is an effective
$1$-cycle,
\begin{multline*}
\label{eq-surf-ocenka-1-new} 0\le -K_X\cdot B\cdot f^*\Sigma=
\\
(2H+f^*(-K_Z-c_1))\cdot (H+f^*(-K_Z-c_1))\cdot f^*\Sigma=
\\
2H^2\cdot f^*\Sigma+ 3 (-K_Z-c_1)\cdot \Sigma=
\\
2 c_1\cdot \Sigma+3 (-K_Z-c_1)\cdot \Sigma= -3K_Z\cdot
\Sigma-c_1\cdot \Sigma.
\end{multline*}
This gives us $0\le a\le b\le 6$.

Moreover, if $a=b=6$, then
\begin{multline*}
-K_X\cdot B\cdot f^*(\Sigma+l)=
\\
(2H+f^*(-K_Z-c_1))\cdot (H+f^*(-K_Z-c_1))\cdot f^*(\Sigma+l)=
\\
2H^2\cdot f^*\Sigma+ 3 (-K_Z-c_1)\cdot (\Sigma+l)=
\\
=2 c_1\cdot (\Sigma+l)-24 =0.
\end{multline*}
Hence, $L\cap B_0=\emptyset$. It follows that $-K_X$ is nef. Indeed,
otherwise there is a curve $R$ such that $K_X\cdot R>0$. Then
$L\cdot R<0$ and $B\cdot R<0$, so $R\subset L\cap B$ and $R\cdot
f^*(\Sigma+l)>0$, a contradiction.
Since $L\cap B_0=\emptyset$, the
map $\Psi\colon X\dashrightarrow V$ given by $\LLL\subset |-K_X|$ is
a morphism near $B_0$ and contracts $B_0$. Moreover, since $V$ has
only isolated singularities the image of $B_0$ is a point. In this
situation the vector bundle $\EEE$ must be decomposable.
This contradicts our assumption.

Thus we may assume that $0\le a\le 5$. Put
\[
p:=\down{a/2}+1,\quad q:=\down{b/2}+1,\quad a':=a-2p,\quad b':=b-2q.
\]
Then
\[
-2\le a',b'\le -1.
\]
Consider the twisted bundle $\EEE':=\EEE\otimes\OOO(-p\Sigma-q l)$.
We have
\begin{equation}
\label{eq-surf-p1-p1-c1-c2} c_1(\EEE')=a'\Sigma+b'l,\quad
c_2(\EEE')=c-aq-bp+2pq.
\end{equation}

\begin{claim}
\label{claim-surf-c-2} We have $c_2(\EEE')<-3$.
\end{claim}
\begin{proof}
Assume that $c_2(\EEE')\ge -3$. Then 
\[
c=c_2(\EEE')+aq+bp-2pq\ge
aq+bp-2pq-3. 
\]
By Riemann-Roch (*) we have
\begin{multline*}
34\le 2g\le \chi(\EEE)\le (a+1)(b+1)-aq-bp+2pq+4
\\
= \frac12 ab +a+b +5+\frac 12 a'b' = \frac 12(a+2)(b+2)+3+\frac 12
a'b'\le 33,
\end{multline*}
a contradiction.
\end{proof}

\begin{claim}
\label{claim-surf-chi} $\chi(\EEE')>0$.
\end{claim}
\begin{proof}
Put $c':=c(\EEE')$. By Riemann-Roch we have
\[
\chi(\EEE')= b'(a'+1)+a'-c'+2.
\]
Assume that $\chi(\EEE')\le 0$. Taking into account Claim
\xref{claim-surf-c-2} we obtain
\[
b'(a'+1)+a' \le c'-2< -5.
\]
Hence $a'\neq -1$. Therefore, $a'=-2$ and
\[
\chi(\EEE')= -b'-c'> 3-b'> 0.
\]
\end{proof}

\begin{claim}
$H^0(\EEE')\neq 0$.
\end{claim}
\begin{proof}
Assume that $H^0(\EEE')=0$. By Claim \xref{claim-surf-chi} we have
$H^2(\EEE')\neq 0$. By Serre Duality
\[
H^2(\EEE')^{*}\simeq H^0(\EEE'^*\otimes \omega_Z)\simeq
H^0(\EEE'\otimes \det \EEE'^*\otimes \omega_Z).
\]
On the other hand, $(\det \EEE'^*\otimes
\omega_Z)^*=\OOO_Z((a'+2)\Sigma+(b'+(e+2))l)$ and $H^0((\det
\EEE'^*\otimes \omega_Z)^*)\neq 0$, a contradiction.
\end{proof}

Now we finish the proof of Proposition \xref{prop-surf-main}.
Consider a nonzero section $s\in H^0(\EEE')$. If $s$ does not vanish
anywhere, then $\EEE'$ is an extension of some line bundle $\EEE_1$
by $\OOO$. But then $c_2(\EEE')=c_2(\EEE_1)=0$. This contradicts
Claim \xref{claim-surf-c-2}. Therefore, the zero locus of $s$
contains a curve $Y$. Let $Y\sim q_1 \Sigma+q_2 l$. Then
restrictions $\EEE'$ to general curves $l\in |l|$ and $l'\in
|\Sigma|$ are $\EEE'|_{l}=\OOO(q_1)\oplus\OOO(a'-q_1)$ and
$\EEE'|_{l'}=\OOO(q_2)\oplus\OOO(b'-q_2)$. Since $-a',-b'\ge 1$, at
least one of the following holds: $2q_1-b'\ge 3$, or $2q_2-a'\ge 3$.
This contradicts Lemma \xref{lemma-p1-p2-nachalo}. Proposition
\xref{prop-surf-main} is proved.

\def\cprime{$'$}


\begin{thebibliography}{10}

\bibitem{Fano-1938-e}
G. Fano.
\newblock {Sulle varieta algebriche a tre dimensioni le cui sezioni iperpiane
  sono superficie di genere zero e bigenere uno.}
\newblock {\em Mem. Mat. Sci. Fis. Natur. Soc. Ital. Sci., III. Ser.},
  24:41--66, 1938.

\bibitem{Godeaux-1933}
L. Godeaux.
\newblock Sur les vari\'et\'es alg\'ebriques \`a trois dimensions dont les
  sections hyperplanes sont des surfaces et de bigenre un.
\newblock {\em Bull. Acad. Belgique Cl. Sci.}, 14:134--140, 1933.


\bibitem{Conte-1983}
A. Conte.
\newblock Two examples of algebraic threefolds whose hyperplane sections are
  {E}nriques surfaces.
\newblock In {\em Algebraic geometry---open problems (Ravello, 1982)}, volume
  997 of {\em Lecture Notes in Math.}, pages 124--130. Springer, Berlin, 1983.

\bibitem{Conte-Murre-1985}
A.~Conte and J.~P. Murre.
\newblock Algebraic varieties of dimension three whose hyperplane sections are
  {E}nriques surfaces.
\newblock {\em Ann. Scuola Norm. Sup. Pisa Cl. Sci. (4)}, 12(1):43--80, 1985.



\bibitem{Bayle-1994}
L. Bayle.
\newblock Classification des vari\'et\'es complexes projectives de dimension
  trois dont une section hyperplane g\'en\'erale est une surface d'{E}nriques.
\newblock {\em J. Reine Angew. Math.}, 449:9--63, 1994.

\bibitem{Sano-1995}
T. Sano.
\newblock On classifications of non-{G}orenstein {${\bf Q}$}-{F}ano {$3$}-folds
  of {F}ano index {$1$}.
\newblock {\em J. Math. Soc. Japan}, 47(2):369--380, 1995.

\bibitem{Minagawa-1999}
T. Minagawa.
\newblock Deformations of {${\bf Q}$}-{C}alabi-{Y}au {$3$}-folds and {${\bf
  Q}$}-{F}ano {$3$}-folds of {F}ano index {$1$}.
\newblock {\em J. Math. Sci. Univ. Tokyo}, 6(2):397--414, 1999.

\bibitem{Cheltsov-1996a}
I.~A. Chel{\cprime}tsov.
\newblock Singularity of three-dimensional manifolds possessing an ample
  effective divisor---a smooth surface of {K}odaira dimension zero.
\newblock {\em Mat. Zametki}, 59(4):618--626, 640, 1996.

\bibitem{Prokhorov-2005a}
Yu.~G. Prokhorov.
\newblock On the degree of {F}ano threefolds with canonical {G}orenstein
  singularities.
\newblock {\em Russian Acad. Sci. Sb. Math.}, 196(1):81--122, 2005.

\bibitem{Utah}
J.~Koll{\'a}r, editor.
\newblock {\em Flips and abundance for algebraic threefolds}.
\newblock Soci\'et\'e Math\'ematique de France, Paris, 1992.
\newblock Papers from the Second Summer Seminar on Algebraic Geometry held at
  the University of Utah, Salt Lake City, Utah, August 1991, Ast\'erisque No.
  211 (1992).

\bibitem{Jahnke-Radloff-2004}
P. Jahnke and I. Radloff.
\newblock Gorenstein Fano threefolds with base points in the anticanonical
  system, 2004.

\bibitem{Przhiyalkovskii-Cheltsov-Shramov-2005}
V.~V. Przhiyalkovski{\u\i}, I.~A. Chel{\cprime}tsov, and K.~A. Shramov.
\newblock Hyperelliptic and trigonal {F}ano threefolds.
\newblock {\em Izv. Ross. Akad. Nauk Ser. Mat.}, 69(2):145--204, 2005.

\bibitem{KM-1998}
J. Koll{\'a}r and S. Mori.
\newblock {\em Birational geometry of algebraic varieties}, volume 134 of {\em
  Cambridge Tracts in Mathematics}.
\newblock Cambridge University Press, Cambridge, 1998.
\newblock With the collaboration of C. H. Clemens and A. Corti, Translated from
  the 1998 Japanese original.

\bibitem{Cossec-Dolgachev-1989-book}
F.~R. Cossec and I.~V. Dolgachev.
\newblock {\em Enriques surfaces. {I}}, volume~76 of {\em Progress in
  Mathematics}.
\newblock Birkh\"auser Boston Inc., Boston, MA, 1989.

\bibitem{Sano-1996}
T. Sano.
\newblock Classification of non-{G}orenstein {${\bf Q}$}-{F}ano {$d$}-folds of
  {F}ano index greater than {$d-2$}.
\newblock {\em Nagoya Math. J.}, 142:133--143, 1996.

\bibitem{Iskovskikh-1979s-e}
V.~A. Iskovskikh.
\newblock Minimal models of rational surfaces over arbitrary fields.
\newblock {\em Math. USSR-Izv.}, 14(1):17--39, 1980.

\end{thebibliography}

\end{document}